# Predictive Control for Energy Management in Ship Power Systems under High-power Ramp Rate Loads

Tuyen V. Vu, *Member, IEEE*, David Gonsoulin, *Student Member, IEEE,* Fernand Diaz, *Student Member, IEEE,* Chris S. Edrington, *Senior, Member, IEEE*, and Touria El-Mezyani, *Member, IEEE*

*Abstract*-- Electrical weapons and combat systems integrated into ships create challenges for their power systems. The main challenge is operation under high-power ramp rate loads, such as rail-guns and radar systems. When operated, these load devices may exceed the ships generators in terms of power ramp rate, which may drive the system to instability. Thus, electric ships require integration of energy storage devices in coordination with the power generators to maintain the power balance between distributed resources and load devices. In order to support the generators by using energy storage systems, an energy management scheme must be deployed to ensure load demand is met. This paper proposes and implements an energy management scheme based on model predictive control to optimize the coordination between the energy storage and the power generators under high-power ramp rate conditions. The simulation and experimental results validate the proposed technique in a reduced scale, notional electric ship power system.

## NOMENCLATURE

| | |
|---|---|
| AES | All-electric ships |
| $A_{ieq}$ | Inequality matrix |
| $b_{ieq}$ | Inequality vector |
| ES | Energy storage device |
| $E_{ES}$ | Energy of energy storage device |
| $E_{ES}^*$ | Energy reference of energy storage device |
| $\bar{E}_{ES}$ | Prediction energy vector of energy storage device |
| $GEN_i$ | Generator $i$ |
| $I_{GENi}$ | Supplying current of generator $i$ |
| $I_{ES}$ | Supplying current of energy storage |
| $I_{PPL}$ | Supplying current of pulsed-power load |
| $I_{Pr}$ | Supplying current of propulsion motor |
| N/A | Not Applicable |
| PPL | Pulsed-power load |
| Pr | Propulsion motor |
| $P_{GENi}$ | Power of generator $i$ |
| $P_{GEN}^{min}$ | Minimum power operation of generators |
| $P_{GEN}^{max}$ | Maximum power operation of generators |
| $P_{ES}$ | Power of energy storage |
| $P_{ES}^{min}$ | Minimum power operation of energy storage |
| $P_{ES}^{max}$ | Maximum power operation of energy storage |
| $\bar{P}_{ES}$ | Prediction power vector of energy storage device |
| $P_{Lj}^*$ | Power reference of load device $j$ |
| $P_{Lj}$ | Power operation of load device $j$ |
| $r_{GENi}$ | Power ramp rate of generator $i$ |
| $r_{GEN}^{min}$ | Minimum power ramp rate of generators |
| $r_{GEN}^{max}$ | Maximum power ramp rate of generators |
| $r_{ES}$ | Power ramp rate of energy storage |
| $r_{ES}^{min}$ | Minimum power ramp rate of energy storage |
| $r_{ES}^{max}$ | Maximum power ramp rate of energy storage |
| $V_{GENi}$ | Terminal voltage at generator $i$ |
| $V_{ES}$ | Terminal voltage at energy storage |
| $V_{GENi}^*$ | Voltage reference at generator $i$ |
| $V_{ES}^*$ | Voltage reference at energy storage |
| $V_{PPL}$ | Terminal voltage at pulsed-power load |
| $V_{Pr}$ | Terminal voltage at propulsion motor |
| $\Delta \bar{P}_{ES}$ | Incremental power vector of energy storage device |
| $\Lambda$ | Lagrange multiplier vector |

This work was supported by the Office of Naval Research via grant awards: N00014-10-1-09.
T. V. Vu is with the Center for Advanced Power Systems – Florida State University, FL 32310 (e-mail: tvu@caps.fsu.edu).
D. Gonsoulin is with the Center for Advanced Power Systems – Florida State University, FL 32310 (e-mail: dg11g@my.fsu.edu).
F. Diaz is with the Center for Advanced Power Systems – Florida State University, FL 32310 (e-mail: fediazfr@gmail.com).
C. S. Edrington is with the Center for Advanced Power Systems – Florida State University, FL 32310 (e-mail: csedrington@gmail.com).
T. El-Mezyani is with the Engineering Department – University of West Florida, FL 32514 (e-mail: telmezyani@uwf.edu).

## I. INTRODUCTION

### A. Literature Review

Future warships will be equipped with advanced systems which integrates electrical weapons and combat systems, such as electromagnetic railguns [1], air and mission defense radars [2], high-energy military lasers [3], and electrical propulsion motors (Pr) [4]. The increase of electric power demand for these devices along with the limited energy capability in the local electric power generations for essential supports has prompted the development of the next generation integrated power systems (NGIPS) for the future naval warships [5]-[7]. In this evolutionary trend, the DC distribution systems for shipboard power systems (SPS) indicated in the NGIPS technology roadmap are under consideration for future warship generations because of the advantages when compared to the AC distribution systems. These competitive advantages include the reduction of power loss and control parameters as well as less complicated power analysis [8].

The DC ship power systems are islanded DC microgrids where the distributed power generators provide the electrical power to their load devices through a DC distribution system.

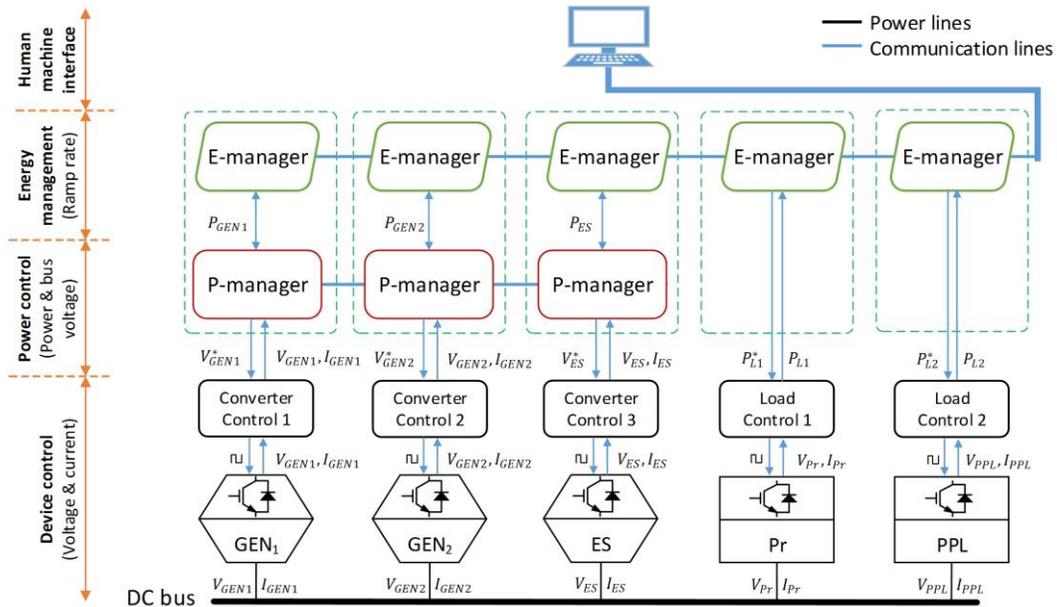

Fig. 1. Distributed control architecture for a notional electric ship.

The main challenge for the power systems is that the essential load devices, such as radars and rail-guns, may have a higher power ramp rate demand than the capability of the power generators. If this condition occurs, it produces unbalanced power between generators and loads, which consequently leads the system to instability. The authors in [9] showed that the generators in real-time simulation can be relaxed in their ramp-rate to handle the constant power load. However, [10] shows a reality that the normal generators such as gas-turbine generators have the ramp rates of 35 to 50 MW/minute, which is far less than the pulsed-power load dynamic operation of 100 MW/second. Therefore, we anticipated that the higher load ramp-rate could drive the system into instability not only because of the voltage drop but also because of the protection mechanism inside the generator units. Hence, the coordination of energy storage (ES) is crucial. When such a high power demand from load devices occurs, the ES will compensate for the power shortage to ensure that the power balance relationship is maintained and the system is stable. The cooperation between generators and ES require an energy management routine to maintain the robustness and stability as well as the optimal operation of the power system.

There are existing energy management strategies for ship power systems [11]-[18]. Among those, the energy management problem for ship power systems has been addressed utilizing a predictive control strategy [12]-[18]. [12] presented an energy management scheme for a system that incorporates two generators. However, ES was not included and no high power ramp rate operations are discussed in the paper. [13] studied the predictive power management in multi-time scale property based on the dynamics of systems. However, the fuel cell used in the work did not focus on supporting the generators in high-power demand cases. Other efforts in [14], [15] presented a predictive control method that considers multiple objectives, such as bus voltage regulation and load demand fulfillment; however, they did not consider the ES as an essential support for synchronous generators in the system. Additionally, these papers considered multiple objectives including slow dynamic terms (power and energy management) and fast dynamic terms (voltage control) for optimization problem. The consideration, however, may result in a computational burden to the central controller because the voltage control requires a small sampling step. Authors in [16] presented the use of ES to mitigate the power fluctuation caused by the Pr; nevertheless, the control algorithm was not designed for an entire ship power system. [17] analyzed the support of ES for pulsed power load operation; yet, the paper did not provide a clear energy management method.

*B. Contribution of the Paper*

Previous predictive control studies in ship power system have not implemented ES to compensate the unbalanced power between generators and high-power ramp rate loads in an experimental system. Therefore, in this paper, we first classify the control system into hierarchical levels, in which we focus on the highest level of control. Second, we consider the ES as an asset to the ship power systems for the control formulation. Finally, we propose a hybrid-MPC, which is unique to solve the fast response requirement for the for the optimized coordination between distributed ES and power generators to satisfy the high-power ramp rate demands of the Pr and the pulsed power load (PPL). Results of the proposed control method are the optimized power references for power generation devices in the ship. The simulation and experimental results for typical study cases of a reduced-scale notional electric ship system are analyzed to validate the effectiveness of the proposed technique.

This paper is organized as follows: Section II introduces the proposed energy management methodology, which includes (1) a hierarchical control architecture and (2) a hybrid predictive energy management algorithm for optimized collaboration between ES and generators. A reduced scale AES with cases is studied in Section III. This Section then provides simulation and experimental results based on these test cases. The results are analyzed and discussed in Section IV to verify the effectiveness of the method. Section V concludes the achievements of the paper.

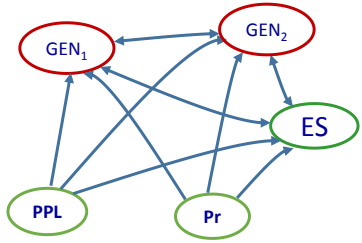

Fig. 2. Communication network among energy managers.

TABLE I
INFORMATION EXCHANGE AMONG ENERGY MANAGERS

| | Send | Receive | Available information |
|---|---|---|---|
| $GEN_1$ | $P_{GEN1}$ | $E_{ES}, P_{ES}, \Delta P_{ES}$ $P_{GEN2}, \Delta P_{PPL}, \Delta P_{Pr}$ | $P_{GEN1}^{min}, r_{GEN1}^{min}$ $P_{GEN1}^{max}, r_{GEN1}^{max}$ |
| $GEN_2$ | $P_{GEN2}$ | $E_{ES}, P_{ES}, \Delta P_{ES}$ $P_{GEN1}, \Delta P_{PPL}, \Delta P_{Pr}$ | $P_{GEN2}^{min}, r_{GEN2}^{min}$ $P_{GEN2}^{max}, r_{GEN2}^{max}$ |
| ES | $E_{ES}, P_{ES},$ $\Delta P_{ES}$ | $P_{GEN1}, P_{GEN2}$ $\Delta P_{PPL}, \Delta P_{Pr}$ | N/A |
| PPL | $\Delta P_{PPL}$ | N/A | N/A |
| Pr | $\Delta P_{Pr}$ | N/A | N/A |

## II. PROPOSED ENERGY MANAGEMENT METHODOLOGY

This section is formatted as follows: First, a control and management architecture for AES is introduced. Second, a hybrid energy management approach based on the architecture is proposed to deal with the power ramp rate problems in AES.

### A. Control and Management Architecture

There are many power and load devices integrated in AES. However, this paper focuses on a notional AES that includes five critical devices, which are a main generator ($GEN_1$), an auxiliary generator ($GEN_2$), an ES, a Pr, and a PPL.

A control and management architecture, which employs control and management schemes for AES is presented first. To manage the operation and coordination of these devices, we configure a three-layer control architecture (Fig. 1). The control architecture consists of the device control, the power control, and the energy management control layers. The first layer includes multiple device controllers, which control the output voltage of generators, output voltage of ES, and the power of load devices (Pr, PPL). The second layer includes multiple power managers (P-managers) to regulate power commands and to maintain bus voltage stability [19]-[21]. The third layer includes energy managers (E-managers) for energy management, which is the focus of this paper. The E-managers employ an energy management scheme to optimize the coordination among the generators and ES to satisfy the load demands.

The energy management scheme for E-managers requires a communication network in which E-managers exchange their power and power ramp rate information. The communication network in this work is configured as shown in Fig. 2, where the Pr and PPL demands are broadcasted to the E-managers of $GEN_1$, $GEN_2$, and the ES. The power and energy information exchanged among energy managers are listed in TABLE I.

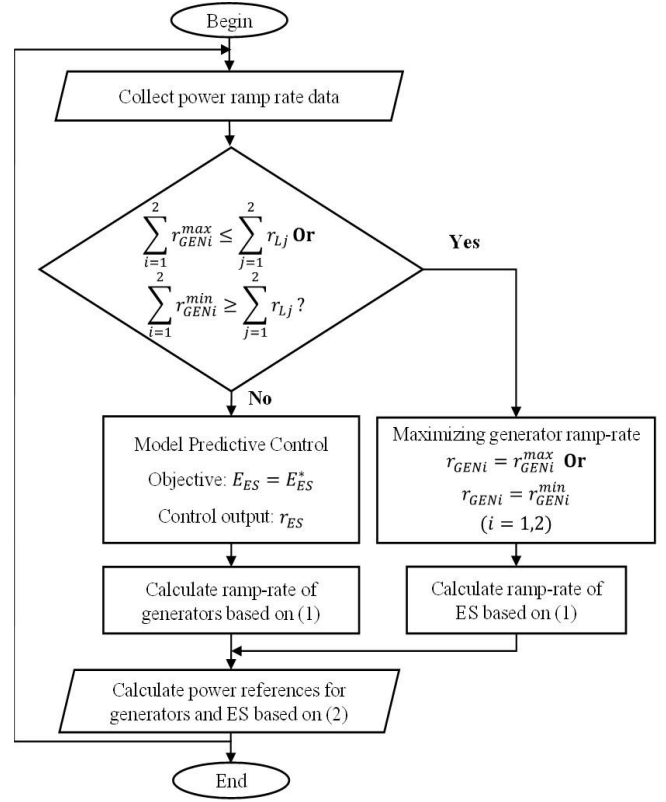

Fig. 3. Hybrid energy management scheme.

### B. Proposed Energy Management Methodology

#### 1) Hybrid energy management approach

The control method is applied in this research as a hybrid approach. The hybrid approach combines the heuristics and MPC (Fig. 3). The heuristics are applied to distinguish system state transition while predictive control is applied fulfil the control objective in each state. In detail, the algorithm first evaluates the aggregated power ramp rate limitation of generators $\sum_{i=1}^{2} r_{GENi}^{min}$, $\sum_{i=1}^{2} r_{GENi}^{max}$ compared to the requested ramp rate of load devices $\sum_{j=1}^{2} r_{Lj}$. If the ramp rate of load devices exceeds the capability of generators, the generators provide their limit power and the ES compensates the unbalanced ramp rate between loads and generators using the equality constraint expressed in (1). If the generators are able to provide ramp rate for load devices, the MPC is utilized to recover the state-of-charge (SOC) of the ES to a predefined reference value $E_{ES}^*$. The output of the predictive method is the optimal ramp rate $r_{ES}$ applied to the ES. Consequently, ramp rate requirements $r_{GEN1}$, $r_{GEN2}$ of generators are calculated via

$$\sum_{j=1}^{2} r_{Lj} = \sum_{i=1}^{2} r_{GENi} + r_{ES}. \quad (1)$$

As described, the hybrid method implemented in E-managers results in the required ramp rate $r_{GEN1}$, $r_{GEN2}$, and $r_{ES}$ applied to the generators and storage devices. Implementation of these ramp rates is conducted through the power commands $P_{GEN1}$, $P_{GEN2}$, and $P_{ES}$, which are calculated based on

$$P_{GENi} = w_{GENi} \int \sum_{i=1}^{2} r_{GENi}$$

$$P_{ES} = \int r_{ES}, \quad (2)$$

where $w_{GENi}$ is the constant parameter proportional to the ramp rate limitation of generator $i$.

$$w_{GENi} = \frac{r_{GENi}^{max}}{\sum_{i=1}^{2} r_{GENi}^{max}} \quad (3)$$

*2) System models*

The development of the control method requires a prediction model of the system for $N_p$ horizon and $N_c$ control steps, which is formulated based on: (a) the ES model, (b) the forecasted load demands, and (c) the predicted power limits in the generators. The control objective is to maintain the SOC of the ES at a desired level and to satisfy the physical constraints in the power system, including the power balance and ramp rate-constraints.

The energy storage model is formulated as

$$E_{ES,k+1} = E_{ES,k} + TP_{ES,k}$$
$$P_{ES,k+1} = P_{ES,k} + Tr_{ES,k}, \quad (4)$$

where subscript $k$ represents a time instant, and $T$ is the sampling time.

Inequality and equality constraints of the systems are described in (5). Specifically, the inequality ramp rate constraint (5a), the inequality power constraint (5b), the inequality energy constraint (5c), the equality power constraint (5d), and the equality ramp rate constraint (5e) are shown as

$$r_{ES}^{min} = r_{GEN}^{min} - r_L \leq r_{ES} \leq r_{ES}^{max} = r_{GEN}^{max} - r_L \quad (5a)$$

$$P_{ES}^{min} = P_{GEN}^{min} - P_L \leq P_{ES} \leq P_{ES}^{max} = P_{GEN}^{max} - P_L \quad (5b)$$

$$0 < E_{ES} < E_{max} \quad (5c)$$

$$P_{GEN} + P_{ES} = P_L \quad (5d)$$

$$r_{GEN} + r_{ES} = r_L, \quad (5e)$$

where $P_{GEN}$ is the total power of generators, $P_L$ is the total power of loads, $r_{GEN}$ is the total ramp rate of generators, and $r_L$ is the total ramp rate of loads. Since the equality constraints (5d) and (5e) are taken care by the heuristics in the hybrid method explained through Fig. 3, only inequality constraints are taken into account for the predictive strategy.

There are constraints in the power ramp rate of ES. Thus, it is useful to utilize the augmented model, which relates the ramp rate of the ES to the control input of the system. Following [22] and [23], the development of such an augmented model is below:

Define

$$\Delta E_{ES,k} = E_{ES,k} - E_{ES,k-1}$$
$$\Delta P_{ES,k} = P_{ES,k} - P_{ES,k-1} = Tr_{ES,k-1}. \quad (6)$$

Substituting $\Delta E_{ES,k}, \Delta P_{ES,k}$ into (4) gives

$$\Delta E_{ES,k+1} = \Delta E_{ES,k} + T\Delta P_{ES,k}. \quad (7)$$

Define $x_{ES,k} = [\Delta E_{ES,k} \quad E_{ES,k}]^T$ as a so-called augmented state. Therefore, combining (4) and (7) yields the following augmented model:

$$x_{ES,k+1} = Ax_{ES,k} + B\Delta P_{ES,k}$$
$$E_{ES,k} = Cx_{ES,k}, \quad (8)$$

where

$$A = \begin{bmatrix} 1 & 0 \\ 1 & 1 \end{bmatrix}, B = \begin{bmatrix} T \\ 0 \end{bmatrix}, C = \begin{bmatrix} 0 & 1 \end{bmatrix}. \quad (9)$$

*3) Predictive control formulation*

From the augmented model expressed in (8), the prediction model for $N_p$ horizon with $N_c$ control steps to describe the energy output of the storage is formulated as

$$\bar{E}_{ES} = Gx_{ES,k} + \Phi \Delta \bar{P}_{ES}, \quad (10)$$

where
$$\bar{E}_{ES} = [E_{ES,k+1}, E_{ES,k+2}, \ldots, E_{ES,k+N_p}],$$
$$\bar{P}_{ES} = [P_{ES,k}, P_{ES,k+2}, \ldots, P_{ES,k+N_p-1}],$$
$$\Delta \bar{P}_{ES} = [\Delta P_{ES,k}, \Delta P_{ES,k+2}, \ldots, \Delta P_{ES,k+N_p-1}],$$
$$\bar{r}_{ES} = [r_{ES,k}, r_{ES,k+2}, \ldots, r_{ES,k+N_p-1}],$$

$$G = \begin{bmatrix} CA \\ CA^2 \\ \vdots \\ CA^{N_p} \end{bmatrix}, \Phi = \begin{bmatrix} CB & 0 & \cdots & 0 \\ CAB & CB & \cdots & 0 \\ \vdots & \vdots & \ddots & \vdots \\ CA^{N_p-1}B & CA^{N_p-2}B & \cdots & CA^{N_p-N_c}B \end{bmatrix}.$$

The objective of the control problem is to manage the energy of the ES to a reference value $E_{ES}^*$ while satisfying the load demand. This objective is represented as minimizing the quadratic cost function

$$J(\Delta \bar{P}_{ES}) = (\bar{E}_{ES}^* - \bar{E}_{ES})^T (\bar{E}_{ES}^* - \bar{E}_{ES}) + \Delta \bar{P}_{ES}^T I_{N_c \times N_c} \Delta \bar{P}_{ES}. \quad (11a)$$

subject to

$$A_{ieq} \Delta \bar{P}_{ES} \leq b_{ieq}, \quad (11b)$$

where

$$\bar{E}_{ES}^* = [E_{ES}^*]_{N_p \times 1}$$
$$A_{ieq} = [A_{ieq}^1 \quad A_{ieq}^2 \quad A_{ieq}^3]^T$$
$$b_{ieq} = [b_{ieq}^1 \quad b_{ieq}^2 \quad b_{ieq}^3]^T. \quad (12)$$

Substituting (10) to (11a) yields

$$J(\Delta \bar{P}_{ES}) = \frac{1}{2} \Delta \bar{P}_{ES}^T M \Delta \bar{P}_{ES} - 2\Phi^T (\bar{E}_{ES}^* - Gx_k) \Delta \bar{P}_{ES}^T + (\bar{E}_{ES}^* - Gx_k)^T (\bar{E}_{ES}^* - Gx_k), \quad (13)$$

where

$$M = 2(\Phi^T \Phi + I_{N_c \times N_c})$$
$$F = -2\Phi^T (\bar{E}_{ES}^* - Gx_k). \quad (14)$$

The elements of $A_{ieq}$ and $b_{ieq}$ are the predicted matrix and vector derived from the inequality constraints expressed in (5a), (5b), and (5c) as

$$A_{ieq}^1 = [-Tri_{N_c} \quad Tri_{N_c}]^T$$
$$A_{ieq}^2 = [-I_{N_c \times N_c} \quad -I_{N_c \times N_c}]^T$$
$$A_{ieq}^3 = [-\Phi \quad \Phi]^T$$
$$b_{ieq}^1 = [(-P_{ES,k}^{min} + P_{ES,k})\bar{1}_{N_c \times 1} \quad (P_{ES,k}^{max} - P_{ES,k})\bar{1}_{N_c \times 1}]^T$$
$$b_{ieq}^2 = T[-r_{ES}^{min} \bar{1}_{N_c \times 1} \quad r_{ES}^{max} \bar{1}_{N_c \times 1}]^T$$
$$b_{ieq}^3 = [-\bar{1}_{N_p \times 1} E_{ES}^{min} + Gx_k \quad \bar{1}_{N_p \times 1} E_{ES}^{max} - Gx_k]^T, \quad (15)$$

where the matrix $Tri_{N_c}$ and vector $\bar{1}$ notations are expressed as

$$Tri_{N_c} = \begin{bmatrix} 1 & 0 & \cdots & 0 \\ 1 & 1 & \cdots & 0 \\ \vdots & \vdots & \ddots & 0 \\ 1 & 1 & \cdots & 1 \end{bmatrix}_{N_c \times N_c}, \bar{\mathbf{1}} = \begin{bmatrix} 1 & 1 & \cdots & 1 \end{bmatrix}^T. \quad (16)$$

The function $J$, representing the objective, is specified as a primal problem, which is equivalent to the problem (17).

$$\max_{\Lambda \geq 0} \min_{\Delta \bar{P}_{ES}} L(\Delta \bar{P}_{ES}, \Lambda) \text{ as}$$

$$L = \frac{1}{2} \Delta \bar{P}_{ES}^T M \Delta \bar{P}_{ES} + \Delta \bar{P}_{ES}^T + \Lambda^T (A_{ieq} \Delta \bar{P}_{ES} - b_{ieq}), \quad (17)$$

where

$$\Lambda = \begin{bmatrix} \lambda_1 & \lambda_2 & \ldots & \lambda_{4N_c + 2N_p} \end{bmatrix}^T. \quad (18)$$

The solution of the minimum value of $L$ will result in the optimal control input $\Delta \bar{P}_{ES}$ applied in the system. The minimization of the cost function is equivalent to the following partial derivative equations:

$$\frac{\partial L}{\partial \Delta \bar{P}_{ES}} = M \Delta \bar{P}_{ES} + F + A_{ieq}^T \Lambda = 0$$

$$\frac{\partial L}{\partial \Lambda} = A_{ieq} \Delta \bar{P}_{ES} - b_{ieq} = 0, \quad (19)$$

which results in

$$\Delta \bar{P}_{ES} = -M^{-1} (F + A_{ieq}^T \Lambda). \quad (20)$$

Thus, the convex optimization problem $L$ is then formulated as a dual problem, which is not computationally intensive and has simpler constraint parameters when compared to the primary problem.

$$\min_{\Lambda \geq 0} \left( \frac{1}{2} \Lambda^T H \Lambda + \Lambda^T K + \frac{1}{2} b_{ieq}^T M^{-1} b_{ieq} \right), \quad (21)$$

where $H = A_{ieq} M^{-1} A_{ieq}^T$, and $K = b_{ieq} + A_{ieq} M^{-1} F$. This optimization problem is solved iteratively utilizing the Hildreth's quadratic programming. $\text{Min}(L(\Delta \bar{P}_{ES}, \Lambda))$ provides the optimized power change in each time step $\Delta \bar{P}_{ES} = T \bar{r}_{ES}$ for ES. The optimized power change in ES results in an optimized power change in generators according to (1). As a result, the power schedule for $P_{ES}$ and $P_{GENi}$ can be calculated through (2).

## III. CASE STUDIES

In this section, a physical testbed representing a notional ship power system and test cases are described. Then, the simulation and experimental results achieved from these test cases will be presented.

### A. Description of System and Test Cases

We illustrate our proposed control approach in a 6 kW experimental testbed as a reduced version (1:10000) of a notional 60 MW ship power system. The reduced-scale testbed is shown in Fig. 4. Because there are no real gas-turbine generators in the laboratory, the two generators are emulated by two power electronic based racks (Rack 1 and Rack 2). The ES is emulated by an electronic battery system (NHR9200), the Pr is represented by the AC load (NHR4600), the PPL is characterized by a DC load (BK Precision). The power and energy capability of generators and ES are indicated in TABLE II. Altough the ES has 8 kW capacity larger than the two generators, we kept this specification for both simulations and experiments for the consistency. The distributed E-managers are implemented on the National Instruments devices (NI myRIO).

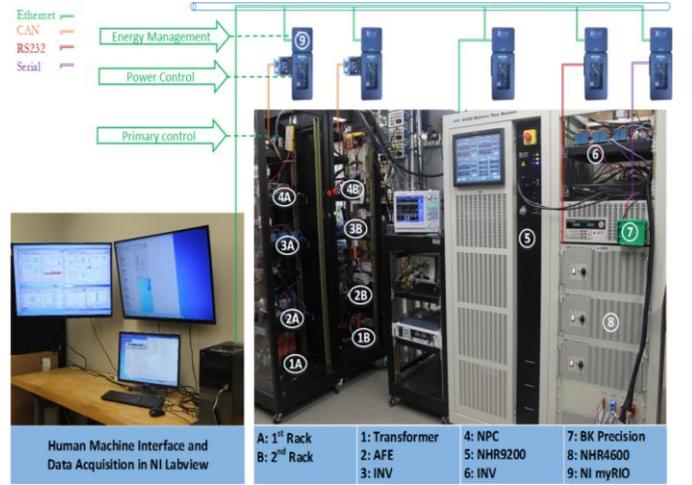

Fig. 4. Reduced-scale AES. Left of the figure shows a Human Machine Interface designed in LabVIEW to acquire information and supervise the myRIO controllers. Right of the figure shows physical testbed with power devices and control devices upper head.

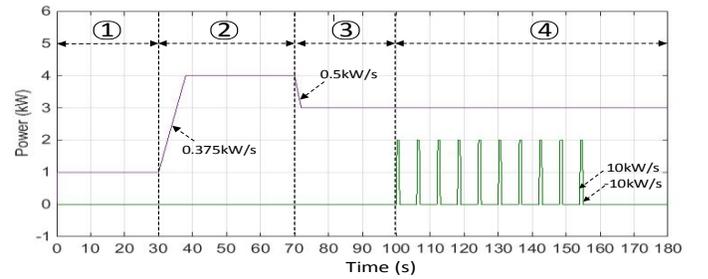

Fig. 5. Load profile.

TABLE II
REDUCED-SCALE AES

| Symbol | Quantity | Values |
|---|---|---|
| $V_B$ | DC bus voltage | 400 V |
| $C_G$ | Capacitor output of a generator | 380 $\mu F$ |
| $P_{GEN1}^{max}$ | Nominal power of GEN$_1$ | 4 kW |
| $P_{GEN2}^{max}$ | Nominal power of GEN$_2$ | 2 kW |
| $r_{GEN1}^{max,min}$ | Power ramp rate limit of GEN$_1$ | $\pm 0.2$ kW/s |
| $r_{GEN2}^{max,min}$ | Power ramp rate limit of GEN$_2$ | $\pm 0.1$ kW/s |
| $E_{ES}$ | Energy storage's capacity | 10 kJ |
| $P_{ES}$ | Energy storage's nominal power | 8 kW |

There are four types of communication setup in the system. The first type is the serial communcation for information exchanging between a myRIO and the BK Precision load underneath. The second type is RS232 for communication between a myRIO and the NHR4600 load underneath. The third type is CAN for communication between a myRIO and an electronic rack underneath. The last type is Ethernet for (a) communication between a myRIO and NHR9200 underneath, and for (b) communication among myRIO devices. The proposed control algorithm is embed with a time step as $T = 10$ ms in the myRIO device using LabVIEW programming. U.S. Navy has specified that the future technology should be ready for 10 rounds per minute firing rate of electromagnetic rail-guns (PPL) [24]. Following that specification, we manage the pulse power load operate every six seconds. Each operation cycle is 1 s. The other 5 s will be devoted for recharging the ES device. To manage the optimal charging process in 5 s with 10 ms time steps we need to look 500 steps ahead, which results in $N_p = 500$. The control horizon we selected as $N_c = 1$ because the control horizon requirement is to be smaller than the prediction horizon. Moreover, $N_c = 1$ will reduce the computation effort for the real-time computation performed by the controller device. The energy, power, voltage, and current data from the testbed are acquired by a Human Machine Interface (HMI) designed using LabVIEW. Moreover, all of the activities of the power devices in the testbed are supervised by that HMI.

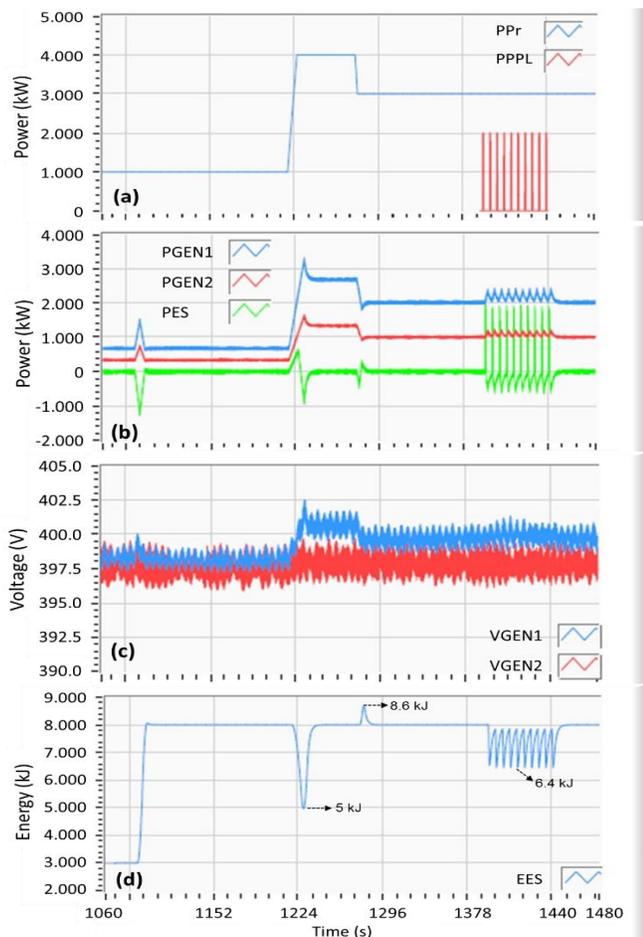

Fig. 6. Data acquisition in NI LabVIEW. (a) load power, (b) generators and ES power, (c) generators terminal voltages, and (d) energy state of ES.

In order to verify the proposed control algorithm, a mission for the notional AES is considered. There are missions, which are described through load profiles (radar, rail-gun, and propulsion) [25], [26]. Similar to these papers, we consider a mission, which consists of four sequential steps for the propulsion and pulsed-power load profiles as shown in Fig. 5. The first step represents a steady state operation of the system, in which the ship moves at constant speed and requires a constant power of 1kW. In this state, the energy of the 10 kJ ES is managed in order to charge from its initial condition at 30% SOC (3 kJ) to a predefined, desired energy level at 80% SOC (8 kJ). The second and third steps reflect an increment of 0.375 kW/s and a decrement of -0.5 kW/s on the propulsion power demand to reach a higher and then lower speed, respectively. In the final step, the ship transfers to a combat mode where it fires the 1.6 kJ rail-gun at the power rate of 10 kW/s for 10 times in 1 minute (note all values are in terms of the scaled test bed).

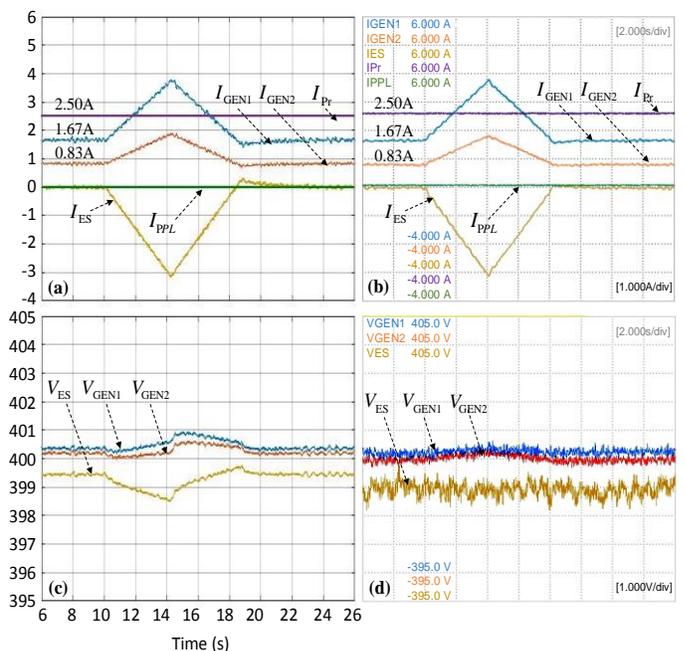

Fig. 7. Voltage and current in stage 1. (a) currents in simulation, (b) experimental currents, (c) voltages in simulation, (d) experimental voltages. [Current (A): 1A/div, Voltage (V): 1V/div, Horizontal axis (s): 2s/div].

*B. Results*

The aforementioned test cases are first implemented in Matlab/Simulink with a virtual communication network with the load profile indicated in Fig. 5. Then the control algorithm is embedded in myRIO devices of ES and generators. In the experiment, there is a load profile embedded in myRIO of PPL and Pr loads; however, there is no restriction about the time for the activation of each load; the activation of each load is remotely controlled through the HMI interface.

The experimental data of the aforementioned test cases for load power, generators and storage power, terminal voltages of generators, and SOC of the ES are acquired by LabVIEW data acquisition, which is presented by the HMI interface as shown in Fig. 6. Details of simulation in Matlab/SIMULINK and experimental validation for the four test cases are shown through Fig. 7 - Fig. 10. These results will be analyzed and discussed in the next section to verify the efficacy of the proposed control algorithm.

## IV. DISCUSSION

The simulation and experimental results acquired in the four aforementioned stages are studied and discussed in this section.

*A. Stage 1*

In this stage, the initial SOC of the ES is 3 kJ, which is 30 %. Fig. 7a shows that before charging process starts, the two generators share the proportional current as $I_{GEN1} = 1.67$ A and $I_{GEN2} = 0.83$ A as they provide 2.5 A (1 kW) to the Pr. While supplying constant power to the Pr, the MPC starts the charging process at 12 s, where there is a current $I_{ES}$ flowing into the battery. The amount of current is provided through $I_{GEN1}$ and $I_{GEN2}$ of two generators. As seen in Fig. 6 the ES is charged to 80 % of SOC, which is the preferred value. At the time the battery is charged, the two generators keep supplying constant power to the propulsion. Experimental results in Fig. 7b validate the simulation results.

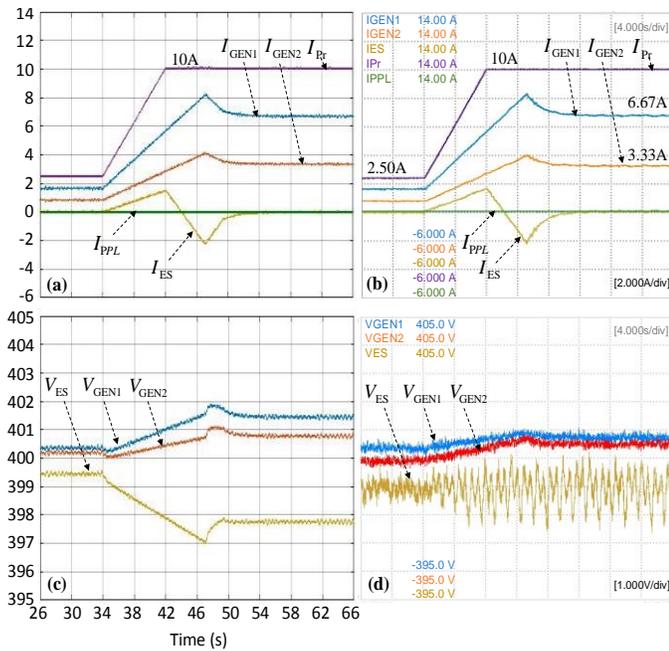

Fig. 8. Voltage and current in stage 2. (a) currents in simulation, (b) experimental currents, (c) voltages in simulation, (d) experimental voltages. [Current (A): 2A/div, Voltage (V): 1V/div, Horizontal axis (s): 4s/div].

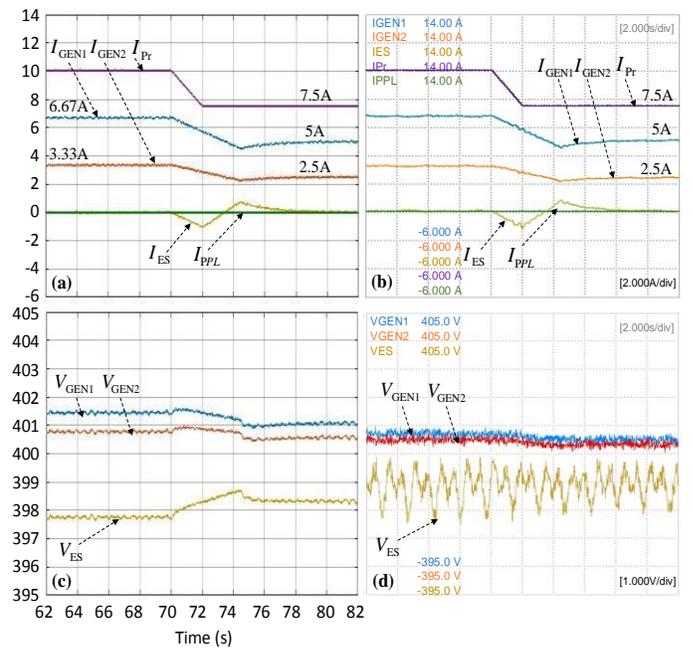

Fig. 9. Voltage and current in stage 3. (a) currents in simulation, (b) experimental currents, (c) voltages in simulation, (d) experimental voltages. [Current (A): 2A/div, Voltage (V): 1V/div, Horizontal axis (s): 2s/div].

Simulation results in Fig. 7c indicate that the terminal voltages $V_{GEN1}$, $V_{GEN2}$, and $V_{ES}$ vary due to the change in current consumption. Specifically, as the currents increase, the voltage deviations from the nominal value at 400 V increase. However, once the charging process is over, the terminal voltages are recovered to the previous values. Experimental results in Fig. 7d shows the terminal voltages of these devices. As seen, these signals contain noise due to the switching and coupling effect of a physical testbed; however, it is shown that the same trend of voltage changes are achieved as the one in simulation. The final voltages are recovered to the initial values.

Consequently, the energy management scheme in this stage is validated through both simulation and experimental results.

### B. Stage 2

This stage continues from the previous stage, and it starts at 34 s in the simulation (Fig. 8). Fig. 8a shows that at 34 s, the ship speeds up as the Pr increases its power consumption from 1 kW to 4 kW at the rate of 0.375 kW/s. The rate of power change is equivalent to the rate of change in the current $I_{Pr}$ as it increases from 2.5 A to 10 A. Due to the fact that this load power ramp rate exceeds the ramp rate capability of the two generators, the ES is set to compensate the unbalanced power ramp rate. As current $I_{Pr}$ reaches 10 A, the MPC recovers the required SOC level of the ES. The ES starts to reduce its power contribution, which also requires the generators to increase their current contribution $I_{GEN1}$ and $I_{GEN2}$. When ES finishes its power support, its energy settles 5 kJ (50 % of SOC) (Fig. 6d). Thus, it requires current supplying from generators for recharging process based on the concurrent constraint values in the system. As the charging process completes, the current of the ES operates at $I_{ES} = 0$ A and the currents of two generators operate at $I_{GEN1} = 6.67$ A and $I_{GEN2} = 3.33$ A. It is also demonstrated through Fig. 8b that the same current values are achieved in the experiment.

Shown in Fig. 8c, the terminal voltage deviations of $V_{GEN1}$, $V_{GEN2}$, and $V_{ES}$ increase as the Pr requires more power. When the predictive method finishes its management, these voltages operate at constant values around 400V as their average value is maintained at 400 V. The voltage deviations are resulted from the increased load condition. In the experiment, Fig. 8d indicates the same trend of terminal voltage deviation. As observed, $V_{ES}$ contains more noise than the terminal voltages of the generators. The noises are from the inverter of the Pr, which is nearby the ES. Therefore, as the power of the Pr increases, the non-linearity of the inverter exhibits noisy response in the DC terminal voltage.

As a result, the energy management strategy for a case of high power ramp rate increment in the Pr has been verified.

### C. Stage 3

In this stage, the Pr reduces its speed as its constant power change decreases from 4 kW(10 A) to 3 kW (7.5 A) with the ramp rate of -0.5 kW/s. This stage starts at 70 s in simulation (Fig. 9). Because the total minimum ramp rate of generators is larger than the load ramp rate, the generators provide their minimum ramp rate. At the same time, the ES supports by absorbing the surplus load power as the current $I_{ES}$ becomes negative to prevent the system from entering an unbalanced power condition (Fig. 9a). The power support is done when the propulsion current $I_{Pr}$ reaches 7.5 A. At this point, the MPC recovers the SOC level of the ES in consideration of ramp rate constraints of generators. As seen, the generators continue to reduce their currents by supplying $I_{GEN1}$ and $I_{GEN2}$ to gradually eliminate ES support until $I_{ES} = 0$ A. As $I_{ES} = 0$ A, the energy level of ES reaches 8.6 kJ (86 % of SOC) (Fig. 6d). The ES continues discharge its energy by increasing its power supplying back to the system. When this discharging process completes, the ES current operates at $I_{ES} = 0$ A and the generator currents are constant values at $I_{GEN1} = 5$ A and $I_{GEN2} = 2.5$ A. The experimental results exhibit the same changing trend and the same final values in ES and generator currents (Fig. 9b).

The terminal voltages $V_{GEN1}$, $V_{GEN2}$, and $V_{ES}$ resulted from simulation are shown in (Fig. 9c). As the propulsion load current decreases, the voltage deviation from 400 V decreases. When the propulsion current reaches the steady state condition and the charging process of the ES is finished, the terminal voltages operate around 400 V. Experimental data in Fig. 9d matches with the simulation results in voltages $V_{GEN1}$ and $V_{GEN2}$ although there is a small mismatch because of the more nonlinearity in the physical system. The voltage $V_{ES}$ seems to have the same behavior as the simulation results; however, the presence of voltage ripple is due to the nearby location of the ES and the Pr.

Hence, simulation and experimental results demonstrated that the control algorithm is effective in cases of a high power decrement in load devices or in cases of large load shedding.

### D. Stage 4

In this last stage, ten pulses are activated with the $\pm 10$ kW/s power ramp rate with the saturated power at 2 kW (5 A). Simulation in Fig. 10a shows that as the PPL increases its power at +10 kW/s, the generators provides its maximum ramp rate capability and the ES provide a large amount of current to compensate the shortage in power ramp rate of the generators. When the PPL saturates its current consumption at $I_{PPL} = 5$ A, the control maximizes the generators' capability to recharge the ES as much as possible. 0.6 s after reaching saturated current value, the PPL changes to the offline mode as it reduces its power to 0 kW at the power rate of -10 kW/s. At this time, the generators reach their ramp rate limits to absorb the large change in ramp rate of PPL, and the ES supports the generators by absorbing the excess power of the PPL. When PPL completes its pulse, the ES current steadies to 0 A, the control method is reactivated again to maximize the recharge process of the ES based on the power ramp rate limit of generators (Fig. 6d). 5 s after the first pulse, the next nine pulses repeat the same power demand as the first pulse. After the PPL finishes ten pulses, the system comes back to the normal operation as generator currents operate at $I_{GEN1} = 5$ A and $I_{GEN2} = 2.5$ A, and the ES current operates at $I_{ES} = 0$ A. The experimental results in Fig. 10b confirms the simulation results as they express the same current behavior and values in every time step.

There are small changes on voltages $V_{GEN1}$, $V_{GEN2}$, and $V_{ES}$ in both simulation (Fig. 10c) and experiment (Fig. 10d). However, this case is expected to exhibit the worse performance of the system. The explanation is that the simulation system and the experimental system utilized the infinite AC sources with power electronics devices to emulate the AC generators. This setup does not truly present the real generators as the power electronics with infinite AC sources have a much faster dynamic than the real generators. However, the current results still verify that the control algorithm optimizes the power scheduling for the generators and ES, which minimize the effect of the large power change in load devices.

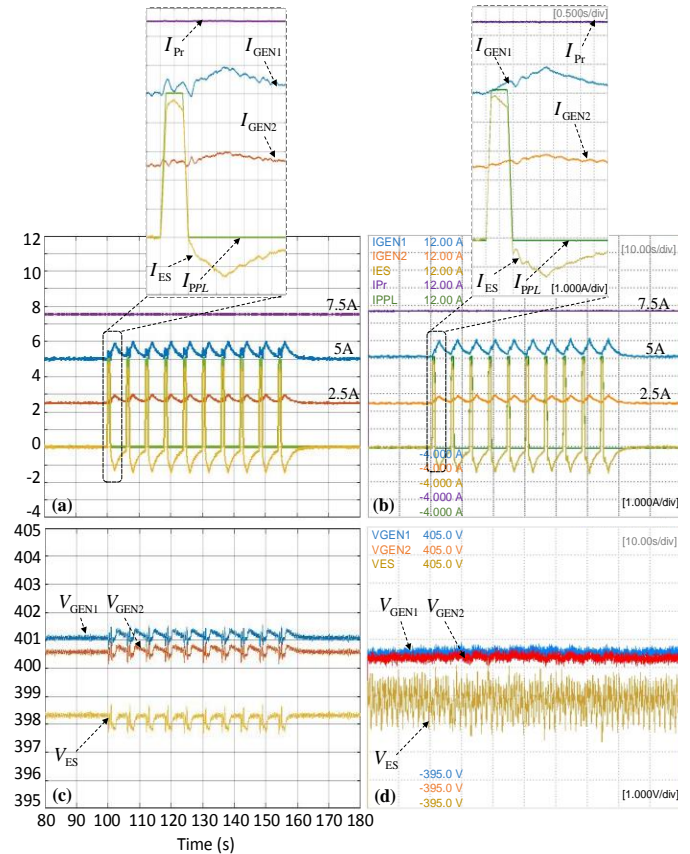

Fig. 10. Voltage and current in stage 4. (a) currents in simulation, (b) experimental currents, (c) voltages in simulation, (d) experimental voltages. [Current (A): 1A/div, Voltage (V): 1V/div, Horizontal axis (s): 10s/div].

## V. CONCLUSION

This paper addresses the operation of AES under high power ramp rate conditions. The problem is solved by a hybrid energy management algorithm, which combines heuristics and MPC. The approach optimizes the use of ES in co-ordination with generators. Simulation and experimental results in a reduced-scaled testbed, which represents a notional AES, have validated the effectiveness of the proposed control algorithm.

In this paper, we assume that ES has enough power ramp-rate to support the operation of pulse power loads. Our future research will consider the limitation in ramp-rates of the ES, which may requires the load-shedding algorithm or the relaxation in the generator units. Furthermore, improving the control algorithm for multiple ES systems will be our interest for our future publications.

## VI. REFERENCES

[1] R. Ellis, "Electromagnetic Railgun," ONR Fact Sheet [Online]. Available:http://www.onr.navy.mil/~/media/Files/Fact%20Sheets/35/Electromagnetic-Railgun-July-2012.ashx
[2] Air and Missile Defense Radar (AMDR), ONR Fact Sheet [Online]. Available:http://www.navy.mil/navydata/fact_display.asp?cid=2100&tid=306&ct=2
[3] R. O'Rourke. (2013, Jan. 22), "Navy Shipboard Lasers for Surface, Air, and Missile Defense: Background and Issues for Congress,"

IEEE Transactions on Energy Conversion 9

Congressional Research Service Report for Congress [Online]. Available: www.fas.org/sgp/crs/weapons/R41526.pdf
[4]  R. O'Rourke. (2013, Feb. 14), "Navy DDG-51 and DDG-1000 Destroyer Programs: Background and Issues for Congress," Congressional Research Service Report for Congress [Online]. Available: http://www.fas.org/sgp/crs/weapons/RL32109.pdf
[5]  Norbert Doerry, " Next Generation Integrated Power Systems (NGIPS) for the Future Fleet," *IEEE Electric Ship Technologies Symposium*, Baltimore, MD April 21, 2009.
[6]  Norbert Doerry, "NGIPS Technology Development Roadmap," NAVSEA Ser 05D / 349 of 30 Nov 2007
[7]  M. Cupelli *et al*., "Power Flow Control and Network Stability in an All-Electric Ship," in *Proceedings of the IEEE*, vol. 103, no. 12, pp. 2355-2380, Dec. 2015.
[8]  J. J. Justo, F. Mwasilu, J. Lee, and J. W. Jung, "AC-microgrids versus DC-microgrids with distributed energy resources: A review*," Renewable and Sustainable Energy Reviews*, vol. 24. pp. 387–405, 2013.
[9]  M. Cupelli, L. Zhu and A. Monti, "Why Ideal Constant Power Loads Are Not the Worst Case Condition From a Control Standpoint," in *IEEE Transactions on Smart Grid*, vol. 6, no. 6, pp. 2596-2606, Nov. 2015.
[10] Lopez, Jaime. "Combustion Engine Vs Gas Turbine- Ramp Rate". *Wartsila.com*. N.p., 2016. Web. 17 Nov. 2016.
[11] K. L. Butler-Purry and N.D.R. Sarma,"Self-healing reconfiguration for restoration of naval shipboard power system," *IEEE Trans. Power Syst.*, vol.19, no. 2, pp. 754-762, May 2004.
[12] S Paran, TV Vu, T El Mezyani, CS Edrington, "MPC-based power management in the shipboard power system," *IEEE Electric Ship Technologies Symposium*, ESTS 2015. pp. 118–122, 2015.
[13] G. Seenumani,"Real-time power management of hybrid power systems in all electric ship applications," *Ph.D. dissertation*, the University of Michigan, 2010.
[14] H. Park, J. Sun, S. Pekarek, P. Stone, D. Opila, R. Meyer, I. Kolmanovsky, and R. DeCarlo, "Real-Time Model Predictive Control for Shipboard Power Management Using the IPA-SQP Approach," *IEEE Transactions on Control Systems Technology*, 2015.
[15] P. Stone, D. F. Opila, H. Park, J. Sun, S. Pekarek ; R. DeCarlo ; E. Westervelt ; J. Brooks ; G. Seenumani, "Shipboard power management using constrained nonlinear model predictive control," *Electric Ship Technologies Symposium (ESTS)*, 2015 IEEE. IEEE, 2015.
[16] J. Hou, J. Sun, and H. Hofmann, "Mitigating Power Fluctuations in Electrical Ship Propulsion Using Model Predictive Control with Hybrid Energy Storage System," *2014 American Control Conference (ACC)*, Portland, Oregon, USA, June 4-6, 2014.
[17] Y. Luo, S. Srivastava, M. Andrus and D. Cartes, "Application of distubance metrics for reducing impacts of energy storage charging in an MVDC based IPS," *2013 IEEE Electric Ship Technologies Symposium (ESTS)*, Arlington, VA, 2013, pp. 287-291.
[18] Y. Luo *et al*., "Application of generalized predictive control for charging super capacitors in microgrid power systems under input constraints," *2015 IEEE International Conference on Cyber Technology in Automation, Control, and Intelligent Systems (CYBER)*, Shenyang, 2015, pp. 1708-1713.
[19] T. V. Vu, S. Paran, F. Diaz Franco, T. El-Mezyani, and C. S. Edrington, "An Alternative Distributed Control Architecture for Improvement in the Transient Response of DC Microgrids," *IEEE Trans. Ind. Electron.*, vol. PP, no. 99, pp. 1–1, Aug. 2016.
[20] T. V Vu, D. Perkins, F. Diaz, D. Gonsoulin, C. S. Edrington, and T. El-Mezyani, "Robust adaptive droop control for DC microgrids," *Electr. Power Syst. Res.*, vol. 146, pp. 95–106, May 2017.
[21] V. Nasirian, S. Moayedi, A. Davoudi, and F. Lewis, "Distributed Cooperative Control of DC Microgrids," *IEEE Trans. Power Electron*., vol. PP, pp. 1–1, 2014.
[22] L. Wang and P.C. Young, "An improved structure for model predictive control using non-minimal state space realization," *Journal of Process Control*, 16 (4):355–371, 2005.
[23] T. V. Vu, S. Paran, F. Diaz, T. E. Meyzani and C. S. Edrington, "Model predictive control for power control in islanded DC microgrids," *Industrial Electronics Society*, IECON 2015 - 41st Annual Conference of the IEEE, Yokohama, 2015, pp. 001610-001615.
[24] "Electromagnetic Railgun Program - Office Of Naval Research". *Onr.navy.mil*. N.p., 2016. Web. 17 Nov. 2016.
[25] A. M. Cramer, X. Liu, Y. Zhang, J. D. Stevens, and E. L. Zivi, ''Early-stage shipboard power system simulation of operational vignettes for dependability assessment,'' in *Proc. 2015 IEEE Electric Ship Technol. Symp. (ESTS),* Alexandria, VA, USA, Jun. 22–24, 2015, pp. 382–387.
[26] D. Stevens, D. F. Opila, A. M. Cramer, and E. L. Zivi, ''Operational vignette-based electric warship load demand,'' *in Proc. IEEE Electr. Ship Technol. Symp.*, Jun. 21–24, 2015, pp. 213–218.

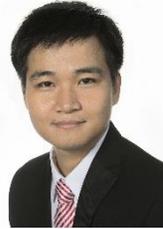

Tuyen V. Vu (S'14, M'17) received his B.S. in electrical engineering from the Hanoi University of Science Technology, Vietnam in 2012, and his Ph.D. in electrical engineering from the Florida State University in 2016.

Since 2016, he has been a postdoctoral research associate in the Florida State University - Center for Advanced Power Systems. His research interests include smart grid; power system dynamics, stability, and control; energy management and optimization; and integration of energy storage systems and electric vehicles into the distribution systems.

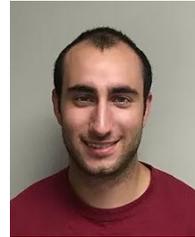

David E. Gonsoulin (S'14) received his BS in engineering from Florida State University in 2015 and is currently pursuing a Ph.D. in Electrical Engineering from Florida State University.

From January 2014 thru May 2015, he was an Undergraduate Research Assistant for the Energy Conversion and Integration Thrust at the Center for Advanced Power Systems. Since 2015, he has been a graduate research assistant at the Florida State University. His research interests include modeling, simulation, and control of terrestrial and ship-based power systems, with a focus on distributed control strategies as well as integration of renewable energy, storage, and pulse power loads.

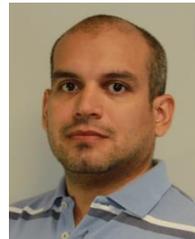

Fernand Diaz Franco (S'14) was born in Cali, Colombia. He received the B.S. degree in physics engineering from University of Cauca, Colombia, and the M.Sc. in mechatronics systems from University of Brasilia, Brazil in 2005 and 2008 respectively.

He is currently working towards the Ph.D. degree in electrical engineering at Florida State University.

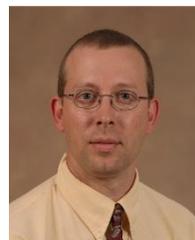

Chris S. Edrington (SM' 08) received his BS in engineering from Arkansas State University in 1999 and his MS and PhD in Electrical Engineering from the Missouri Science and Technology (formerly University of Missouri-Rolla) in 2001 and 2004, respectively, where he was both a DoE GAANN and NSF IGERT Fellow.

He currently is a Professor of Electrical and Computer Engineering with the FAMU-FSU College of Engineering and is the lead for the Energy Conversion and Integration thrust for the Florida State University-Center for Advanced Power Systems. His current research is in the application of distributed controls and optimization for power and energy management in maritime and terrestrial power systems as well as real-time stability and complexity assessment of such highly nonlinear systems.

Dr. Edrington has published over 150 papers (including 2 IEEE Prize Awards), has graduated 25 MS students and 8 Ph.D. students His research has also resulted in 4 patents awarded, with 11 additional patents to be awarded over the next year, as well as 1 patent pending.

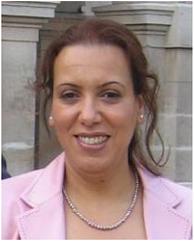 Touria El-Mezyani (M'08) received the Ph.D. degree in Control, Computer and Systems Engineering from Université des Sciences et Technologies de Lille, France in 2005.

From 2013 to 2016, she was a Research Faculty with the Center for Advanced Power Systems, Florida State University, Tallahassee. Since 2016, she has been an assistant professor at the University of West Florida. Her research interests include complex and hybrid systems and their application to power systems, distributed/decentralized controls and optimization, fault detection, isolation, and diagnosis, control and optimization, multi-agent systems. She is an affiliate of the IEEE Automatic Control Society since 2008.